\documentclass[11pt,a4paper]{amsart}
\usepackage{amsfonts}
\usepackage{mathrsfs}
\usepackage{CJK}

\usepackage{palatino}
\usepackage{amsmath,amssymb,amsfonts}
\usepackage[breaklinks,bookmarksopen,bookmarksnumbered]{hyperref}

\allowdisplaybreaks

\theoremstyle{plain}
\newtheorem{thm}{Theorem}[section]
\newtheorem{lem}{Lemma}[section]

\newtheorem{cor}{Corolary}[section]

\numberwithin{equation}{section}

\theoremstyle{definition}

\def\pr{\noindent\textit{Proof}:\quad}

\begin{document}
\title{Rotation invariant singular K\"ahler metrics with constant scalar curvature on $\mathbb{C}^n$}
\author{Weiyong He and Jun Li}
\address{Department of Mathematics, University of Oregon, Eugene, OR 97403. }
\email{whe@uoregon.edu;}
\address{College of Mathematics and Econometrics, Hunan University, Changsha, 410082, P. R. China}
\email{lijun1985@hnu.edu.cn}
\date{\today}
\maketitle

\section{Introduction}

 In this paper we study rotation invariant K\"ahler metric with constant scalar curvature (csck for short) on $\mathbb{C}^n \backslash \{0\}$. One of our main motivations is to understand the singularity of solutions of the scalar curvature type equation in view of existence of Calabi's extremal metric \cite{Ca1}. Equation of this type is very hard and was not understood well in the literature.
           
 To illustrate the main points, we consider a viscosity solution $u(z_1, z_2)=2(1+|z_1|^2)|z_2|$ \cite{He10} of the equation $\det(u_{i\bar j})=1$ on $\mathbb{C}^2$. It is a continuous potential of a singular K\"ahler metric with both zero Ricci curvature and scalar curvature on $\mathbb{C}^2$. The singular locus $\{z_2=0\}$ cannot be trapped inside a compact subset. 
          
 One may speculate that this is a common phenomenon for csck. The csck cannot have  compact singularities if the potential is bounded. In fact if such an example were to exist, one would be able to construct a singular K\"ahler metric with both bounded scalar curvature and potential.

Here we consider a much more modest question: whether a csck on $\mathbb{C}^n$ can have isolated singularity and what is the behavior of the metric near this singularity if it exists? To make this question accessible, we assume moreover that the metric is invariant under rotations.    
           
Our results regarding rotation invariant csck on $\mathbb{C}^n \backslash \{0\}$ confirm the above speculation. The potential and the entropy of such kind metrics cannot be bounded near the singularity $0 \in \mathbb{C}^n$. Indeed all such singular potentials have logarithm poles. 
           
 During the preparation of this paper, Chen and Cheng \cite{CC1, CC2, CC3} have obtained substantial progress in understanding the regularity of scalar curvature type equation on compact K\"ahler manifolds. They proved, for a csck, if its potential is bounded, then all higher derivatives are priori bounded by its potential bound. They also proved that if scalar curvature and the entropy are bounded, the metrics satisfying the bounds are precompact in $C^{1, \alpha}$. With these deep estimates, Chen-Cheng \cite{CC3} solved several substantial conjectures in the field regarding the existence of csck, including the properness conjecture and the geodesic stability conjecture \cite{Don97, Chen06}. Their great achievements certainly answered the problems which motivates us to consider rotation invariant csck. It should be mentioned that the method in \cite{CC1, CC2, CC3} works essentially only for compact manifolds. However we hope our findings are still interesting in its own right.

 Now we describe our results. For a K\"ahler metric $\omega$ on $\mathbb{C}^n\backslash\{0\}$, the scalar curvature type equation is a fourth order PDE on the potential. This PDE can be reduced to a system of ODE \ref{main equation} if the potential is radial. By carefully examining its structure, we can transform the system of ODE into a single ODE \ref{pequation01} if the scalar curvature is a constant. This ODE can be solved completely in the smooth case and lower dimensional case. 
    
  By solving the ODE \ref{pequation01} in the smooth case, we can give a complete list of rotation invariant csck on $\mathbb{C}^n$. Though this result may be well known, we include it here for comparison.
  
    \begin{thm}\label{n smooth solution}
    Suppose $n \geq 2$ is an integer.
\begin{itemize}
\item[(1)] The rotation invariant K\"ahler metric $\omega$ with zero constant scalar curvature on $\mathbb{C}^n$ must be a multiple of the standard Euclidean metric.

\item[(2)] The rotation invariant K\"ahler metric $\omega$ with constant scalar curvature $n(n+1)$ on $\mathbb{C}^n$ must be of the form
\[
		\omega=i\frac{\sum_{j=1}^{n}dz_j \wedge d \overline{z}_j}{\sum_{j=1}^{n}|z_j|^2+a}-i\frac{(\sum_{j=1}^{n}\overline{z}_jdz_j)\wedge \overline{(\sum_{j=1}^{n}\overline{z}_jdz_j)}}{(\sum_{j=1}^{n}|z_j|^2+a)^2}
\]
where $a>0$ is a constant.

\item[(3)] There does not exist rotation invariant K\"ahler metric with negative constant scalar curvature on $\mathbb{C}^n$.
\end{itemize}
\end{thm}

In the positive scalar curvature case, the series of smooth K\"ahler metrics degenerate to a metric with isolated singularity at $0 \in \mathbb{C}^n$ as $a\rightarrow 0$. Along this degeneration, the potential and the entropy both diverge. 
        
By solving the ODE \ref{pequation01} for dimenson $n<4$, we obtain a complete list of radial K\"ahler potentials with zero or positive constant scalar curvature on $\mathbb{C}^2\backslash\{0\}$ and $\mathbb{C}^3\backslash\{0\}$.

\begin{thm}\label{2 zero singular solution}
Let $u:(0,+\infty) \rightarrow \mathbb{R}$ be a smooth function such that $u(|z_1|^2+|z_2|^2)$ is the potential of a K\"ahler metric with constant scalar curvature $R=0$ on $\mathbb{C}^2 \backslash \{0\}$. Then one of the following is true:

\begin{itemize}
\item[(1)] There exist constants $a,b$ with $a>0$ such that
 \[			
 u(s)=as+b
 \]
 		
\item[(2)] There exist constants $a,b,c$ with $a>0 ,b> 0$ such that
\[
 		u(s)=as+b\log s+c
\]
 		
\item[(3)] There exist constants $\alpha,\beta,c$ with $\alpha \neq 0, \beta>0, \alpha < \beta$
 such that $g(s)=su'(s)$ is the smooth strictly increasing function on $(0,+\infty)$ ranging from $\beta$ to $+\infty$   determined by
\[
 			\frac{\beta}{\beta-\alpha} \log(g(s)-\beta) -\frac{\alpha}{\beta-\alpha} \log(g(s)-\alpha)=\log s +c
\]
 		 		
\item[(4)] There exist constants $\alpha,c$ with $\alpha >0$ such that $g(s)=su'(s)$ is the smooth strictly increasing function ranging from $\alpha$ to $+\infty$ on $(0,+\infty)$ determined by
\[
 	   	\log(g(s)-\alpha)-\frac{\alpha}{g(s)-\alpha}=\log s +c
\]
\end{itemize}
\end{thm}
 
\begin{thm}\label{2 negative singular solution}
Let $u:(0,+\infty) \rightarrow \mathbb{R}$ be a smooth function such that $u(|z_1|^2+|z_2|^2)$ is the potential of a K\"ahler metric with constant scalar curvature $R=6$ on $\mathbb{C}^2\backslash\{0\}$ and $g(s)=su'(s)$. Then one of the following is true:

\begin{itemize}
\item[(1)] There exist constants $a,c$ with $a>0$ such that
\[
			u(s)=\log(s+a)+c
\]

\item[(2)] There exist constants $a, k$ with $a>0, 0<k <1$ such that
\[
			g(s) =\frac{1}{2}(k+1)-\frac{k a}{s^k+a}
\]
		
\item[(3)] There exist constants $\alpha,\beta,\gamma,c$ with 
           $\alpha \neq 0, \beta>0, \alpha<\beta<\gamma, \alpha+\beta+\gamma=1$ 
such that $g(s)=su'(s)$ is the smooth strictly increasing function ranging from $\beta$ to $\gamma$ on $(0,+\infty)$ determined by

\begin{equation*}
\begin{split}
    	&-\alpha(\gamma-\beta)\log(g(s)-\alpha)+\beta(\gamma-\alpha)\log(g(s)-\beta) \\
     	&-\gamma(\beta-\alpha)\log(\gamma-g(s)) 
    	=(\beta-\alpha)(\gamma-\beta)(\gamma-\alpha)\log s+c.
 \end{split}	
 \end{equation*}	
		
\item[(4)] There exist  constants $\alpha,\beta,\gamma$ with $0<\beta<\alpha,\alpha+2\beta=1$
such that $g(s)=su'(s)$ is the smooth strictly increasing function ranging from $\beta$ to $\alpha$  on $(0,+\infty)$ determined by
\[
		\alpha\log(g(s)-\beta)-\alpha\log(\alpha-g(s))+\frac{\beta(\beta-\alpha)}{g(s)-\beta}=(\beta-\alpha)^2\log s+c
\]
\end{itemize}
\end{thm}

\begin{thm}\label{3 zero singular solution}
	Let $u:(0,+\infty) \rightarrow \mathbb{R}$ be a smooth function such that $u(|z_1|^2+|z_2|^2+|z_3|^2)$ is the potential of a K\"ahler metric $\omega$ with constant scalar curvature $R=0$ on $\mathbb{C}^3\backslash\{0\}$ and $g(s)=su'(s)$. Then one of the following is true:
		
\begin{itemize}
\item[(1)] There exist constants $a, b$ with $a>0$ such that
\[
			u(s)=as+b
\]

\item[(2)] There exist constants $a, b, c$ with $a, b>0$ such that
\[
			u(s)=a(\sqrt{s^2+b^2}+b\log\frac{s}{\sqrt{s^2+b^2}+b})+c
\]
		
\item[(3)] There exist constants $\alpha,\beta,\gamma$ with 
\[
		\alpha<0, \beta \neq 0, \gamma>0, \alpha<\beta<\gamma, \alpha+\beta+\gamma=0
\] 
such that $g(s)$ is the  smooth strictly increasing function ranging from $\gamma$ to $+\infty$ on $(0,+\infty)$ determined by
\begin{equation*}
\begin{split}
	             	&\alpha^2(\gamma-\beta)\log(g(s)-\alpha)-\beta^2(\gamma-\alpha)\log(g(s)-\beta)  \\
		&+\gamma^2(\beta-\alpha)\log(g(s)-\gamma) 
		=(\beta-\alpha)(\gamma-\beta)(\gamma-\alpha)\log s+c.
\end{split}
\end{equation*}	
	
\item[(4)] There exist constants $\alpha, c$ with $\alpha<0$ such that $g(s)$ is the smooth strictly increasing function ranging from $-2\alpha$ to $+\infty$ on $(0,+\infty)$ determined by
\[
		\frac{5}{9}\log (g(s)-\alpha)+\frac{4}{9}\log(g(s)+2\alpha)-\frac{\alpha}{3}\frac{1}{g(s)-\alpha}=\log s+c
\]
		
\item[(5)] There exist constants $\alpha, c$ with $\alpha>0$ such that $g(s)$ is the smooth strictly increasing function ranging from $\alpha$ to $+\infty$ on $(0,+\infty)$ determined by
\[
			\frac{5}{9}\log(g(s)-\alpha)+\frac{4}{9}\log (g(s)+2\alpha)-\frac{\alpha}{3}\frac{1}{g(s)-\alpha}=\log s+c
\]
		
\item[(6)] There exist constants $\alpha,\beta,\gamma, c$ with 
\[
		\alpha>0,\gamma>0,\alpha+2\beta=0
\] 
such that $g(s)$ is the smooth strictly increasing function ranging from $\alpha$ to $+\infty$ on $(0,+\infty)$ determined by
\begin{equation*}
\begin{split}
                 &\alpha^2\log(g(s)-\alpha)  
                 +\frac{(\beta^2+\gamma^2-2\alpha\beta)}{2}\log[(g(s)-\beta)^2+\gamma^2] \\ 
                 &+\frac{(\beta^3+\beta\gamma^2-\alpha\beta^2+\alpha\gamma^2)}{\gamma}\arctan \frac{g(s)-\beta}{\gamma}=[(\alpha-\beta)^2+\gamma^2]\log s+c
\end{split}	
\end{equation*}

\end{itemize}
\end{thm}

\begin{thm}\label{3 negative singular solution}
Let $u:(0,+\infty) \rightarrow \mathbb{R}$ be a smooth function such that $u(|z_1|^2+|z_2|^2+|z_3|^2)$ is the potential of a K\"ahler metric with constant scalar curvature $R=12$ on $\mathbb{C}^3\backslash\{0\}$ and  $g(s)=su'(s)$. Then one of the following is true:

\begin{itemize}
\item[(1)] There exist constants $a, b$ with $a>0$ such that
\[
		u(s)=\log (s+a)+b
\]
	
\item[(2)] There exist constants $\alpha,\beta,\gamma, c$ with 
\[
		\alpha<0<\beta<\gamma,\alpha+\beta+\gamma=1,\alpha\beta+\beta\gamma+\gamma\alpha=0
\]
such that $g(s)$ is the smooth strictly increasing function ranging from $\beta $ to $\gamma$ on $(0,+\infty)$ determined by
\begin{equation*}
\begin{split}
	       &-\alpha(\gamma-\beta)\log(g(s)-\alpha)+\beta(\gamma-\alpha)\log(g(s)-\beta) \\
	      & -\gamma(\beta-\alpha)\log(\gamma-g(s))=(\beta-\alpha)(\gamma-\beta)(\gamma-\alpha)\log s+c
\end{split}
\end{equation*}

\item[(3)] There exist constants $\alpha,\beta,\gamma,\delta, c$ with
\[
		\alpha<\beta<\gamma<\delta,\alpha+\beta+\gamma+\delta=1,\alpha\beta+\alpha\gamma+\alpha\delta+\beta\gamma+\beta\delta+\gamma\delta=0,\gamma>0,\alpha\beta\delta \neq 0
\]
	 such that $g(s)$ is the smooth strictly increasing function ranging from $\gamma$ to $\delta$ on $(0,+\infty)$ determined by
\begin{equation*}
\begin{split}
	&\frac{\alpha^2}{(\beta-\alpha)(\gamma-\alpha)(\delta-\alpha)}\log(g(s)-\alpha)-\frac{\beta^2}{(\beta-\alpha)(\gamma-\beta)(\delta-\beta)}\log(g(s)-\beta)   \\
	&+\frac{\gamma^2}{(\gamma-\alpha)(\gamma-\beta)(\delta-\gamma)}\log(g(s)-\gamma)-\frac{\delta^2}{(\delta-\alpha)(\delta-\beta)(\delta-\gamma)}\log(\delta-g(s)) \\
	&=\log s+c
\end{split}
\end{equation*}
	
\item[(4)] There exist constants $\alpha,\beta,\gamma, c$ with
\[
		2\alpha+\beta+\gamma=1,\alpha^2+2\alpha\beta+2\alpha\gamma+\beta \gamma ,\alpha<0<\beta<\gamma
\] 
	such that $g(s)$ is the smooth strictly increasing function ranging from $\beta$ to $\gamma$ on $(0,+\infty)$ determined by
\begin{equation*}
\begin{split}
	 &-\frac{\gamma^2}{(\gamma-\beta)(\gamma-\alpha)^2}\log(\gamma-g(s))+\frac{\beta^2}{(\gamma-\beta)(\beta-\alpha)^2}\log(g(s)-\beta) \\
	 &-\frac{2\alpha\beta\gamma-\alpha^2(\beta+\gamma)}{(\gamma-\alpha)^2(\beta-\alpha)^2} \log(g(s)-\alpha)  
	 +\frac{\alpha^2}{(\gamma-\alpha)(\beta-\alpha)}\frac{1}{g(s)-\alpha}=\log s+c
\end{split}
\end{equation*}
	 
\item[(5)] There exist constants $\alpha,\beta,\gamma,c$ with
\[
	 	2\alpha+\beta+\gamma=1,\alpha^2+2\alpha\beta+2\alpha\gamma+\beta\gamma=0,\beta<0<\alpha<\gamma
\]
 such that $g(s)$ is the smooth strictly increasing function ranging from $\alpha$ to $\gamma$ on $(0,+\infty)$ determined by
\begin{equation*}
\begin{split}
	 &-\frac{\gamma^2}{(\gamma-\beta)(\gamma-\alpha)^2}\log(\gamma-g(s))+\frac{\beta^2}{(\gamma-\beta)(\beta-\alpha)^2}\log(g(s)-\beta) \\
	 &-\frac{2\alpha\beta\gamma-\alpha^2(\beta+\gamma)}{(\gamma-\alpha)^2(\beta-\alpha)^2} \log(g(s)-\alpha)  
	 +\frac{\alpha^2}{(\gamma-\alpha)(\beta-\alpha)}\frac{1}{g(s)-\alpha}=\log s+c
\end{split}
\end{equation*}

\item[(6)]  There exist constants $\alpha_1,\alpha_2,\beta,\gamma,c$ with 
\[
	 	0<\alpha_1 <\alpha_2,\gamma>0,\alpha_1+\alpha_2+2\beta=1,\alpha_1\alpha_2+2\beta(\alpha_1+\alpha_2)+\beta^2+\gamma^2=0
\]
	 such that $g(s)$ is the smooth strictly increasing function ranging from $\alpha_1$ to $\alpha_2$ on $(0,+\infty)$ determined by
	 \begin{equation*}
	 \begin{split}
	 &\frac{1}{(\alpha_1-\beta)^2+\gamma^2}[\alpha^2_1\log(g(s)-\alpha_1) +\frac{\beta^2+\gamma^2-2\alpha_1\beta}{2}\log((g(s)-\beta)^2+\gamma^2)  \\ 
	 &+\frac{(\alpha_1+\beta)(\beta^2+\gamma^2)-2\alpha_1\beta^2}{\gamma}\arctan\frac{g(s)-\beta}{\gamma}]  \\
	 &-\frac{1}{(\alpha_2-\beta)^2+\gamma^2}[\alpha^2_2\log(\alpha_2-g(s))+\frac{(\beta^2+\gamma^2-2\alpha_2\beta)}{2}\log((g(s)-\beta)^2+\gamma^2) \\
	 &+\frac{(\alpha_2+\beta)(\beta^2+\gamma^2)-2\alpha_2\beta^2}{\gamma}\arctan\frac{g(s)-\beta}{\gamma}]  \\
	 &=(\alpha_2-\alpha_1)\log s+c
	 \end{split}
	 \end{equation*}
\end{itemize}
\end{thm}

 Potentials of the first case in Theorem \ref{2 zero singular solution}, \ref{2 negative singular solution}, \ref{3 zero singular solution}, \ref{3 negative singular solution} correspond to the smooth K\"ahler metric we obtained in Theorem \ref{n smooth solution}(Lemma \ref{smooth}). Potentials of other cases have a logarithmic pole and hence have nonzero Lelong number at the origin. In particular, the second case in Theorem \ref{2 zero singular solution} corresponds to the Burn's metric \cite{LeB 91}.
 
If the scalar curvature is zero or positive, One can check that all local solutions for \ref{pequation01} near $0$ can be extended to a solution on $(0,+\infty)$ and hence determines a metric on $\mathbb{C}^n \backslash \{0\}$. But if the scalar curvature is negative, every local solution near $0$ can only be extended to a finite interval $(0,s_0)$. In fact we can prove that

\begin{thm}\label{positive}
There does not exist rotation invariant K\"ahler metric with  negative constant scalar curvature  on $\mathbb{C}^2\backslash\{0\}$ or $\mathbb{C}^3\backslash\{0\}$.
\end{thm}

Those local solution near $0$ can only be extended to a finite interval $(0,s_0)$ will lead to singular radial K\"ahler metric with negative constant scalar curvature on the ball $B(0,\sqrt{s_0})$ which can be normalized to a metric on the unit ball. In dimension $2$ we have the following:

\begin{thm}\label{ball}
Suppose $u:(0,1) \rightarrow \mathbb{R}$ be a smooth function such that $u(|z_1|^2+|z_2|^2)$ is the potential of a singular radial K\"ahler metric with constant scalar curvature $R=-6$ on the unit ball $B(0,1) \subset \mathbb{C}^2$ and $g(s)=su'(s)$. Moreover if $g(s)$ satisfies the normalized condition $\lim_{s \rightarrow 1} g(s)=+\infty$, then one of the following is true:
\begin{itemize}
	\item[(1)] $g(s)=\frac{s}{1-s}$.
	\item[(2)] There exists a constant $k>1$ such that
	\[
	g(s)=-\frac{1}{2}(k+1)+\frac{k}{1-s^k}.
	\]
	\item[(3)] There exist constants $\alpha,\beta,\gamma$ with 
	\[
	\alpha+\beta+\gamma=-1,\alpha<\beta<\gamma,\gamma>0
	\]
	such that $g(s):(0,1) \rightarrow (\gamma,+\infty)$ is the strictly increasing function defined on $(0,1)$, ranging from $\gamma$ to $\infty$ determined by
	\[
	(g(s)-\alpha)^{\alpha(\gamma-\beta)}(g(s)-\beta)^{-\beta(\gamma-\alpha)}(g(s)-\gamma)^{\gamma(\beta-\alpha)}=s^{(\beta-\alpha)(\gamma-\alpha)(\gamma-\beta)}.
	\]
	\item[(4)] There exist constants $\alpha,\beta$ such that
	\[
	\alpha<0,\beta>0,\alpha+2\beta=-1
	\]
	and $g(s)$ is the strictly increasing function defined on $(0,1)$, ranging from $\beta$ to $\infty$ determined by
	\[
	-\alpha \log (g(s)-\beta)+\alpha \log(g(s)-\alpha)-\frac{\beta(\beta-\alpha)}{g(s)-\beta}=(\beta-\alpha)^2\log s.
	\]
	\item[(5)] There exist constants $\alpha,\beta$ such that
	\[
	\alpha>0,\beta<0,\alpha+2\beta=-1
	\]
	and $g(s)$ is the strictly increasing function defined on $(0,1)$, ranging from $\alpha$ to $\infty$ determined by
	\[
	-\alpha \log (g(s)-\beta)+\alpha \log(g(s)-\alpha)-\frac{\beta(\beta-\alpha)}{g(s)-\beta}=(\beta-\alpha)^2\log s.
	\]
	\item[(6)]There exist constants $\alpha,\beta,\gamma$ such that
	\[
	\alpha>0,\gamma>0,\alpha+2\beta=-1
	\]
	and $g(s)$ is the strictly increasing function defined on $(0,1)$, ranging from $\alpha$ to $\infty$ determined by
	\[
	\log(g(s)-\alpha)-\frac{1}{2}\log[(g(s)-\beta)^2+\gamma^2]-\frac{\alpha-\beta}{\gamma}\arctan \frac{g(s)-\beta}{\gamma}=[(\alpha-\beta)^2+\gamma^2]\log s-\frac{\pi(\alpha-\beta)}{2\gamma}.
	\]
\end{itemize}
\end{thm}

Theorem \ref{ball} is obtained by a careful check of all local solutions for \ref{pequation01}, but the argument to show that we have exhausted all possiblilities seems to be too tedious. So we omit the proof here. We want to remark that the first case corresponds to the well known Bergman metric on the unit ball.

We believe that Theorem \ref{positive} holds for all dimension. But it is beyond the method we adopted here. In fact the ODE \ref{pequation01} becomes very complicated for dimension $n \geq 4$. It also explains why we fail to give a list of rotation invariant csck on $\mathbb{C}^n \backslash \{0\}$ in zero and positive curvature case for all dimensions. It seems that the solutions may be quite wild.\\

  {\bf Acknowledgement:}  The first named author is supported in part by an NSF grant, award no. 1611797. The second named author wants to thank Prof. Xiangyu Zhou and Prof. Yueping Jiang for help and encouragements. He is partially supported by NSFC 11701164.
         
\section{The constant scalar curvature equation on radial K\"ahler potential}
    
    Let $u:(0,+\infty) \rightarrow \mathbb{R}$ be a smooth function. The real $(1,1)$ form
    \begin{equation*}
    \omega=i\partial \overline{\partial}u(|z|^2)=i\sum_{j,k=1}^{n}(\delta_{jk}u'(s)+u''(s)\overline{z}_jz_k)dz^j\wedge d\overline{z}^k
    \end{equation*}
   is a K\"ahler metric on $\mathbb{C}^n \backslash \{0\}$ if and only if
    \begin{equation*}
    	u'(s)>0,u'(s)+su''(s)>0.
    \end{equation*}
    
    Denote by $G=(g_{j\overline{k}})$ the metric matrix of the K\"ahler metric. Recall that $\det(I_p-MN)=\det(I_q-NM)$ for $p \times q$ matrix M and $q \times p$ matrix $N$ (\cite{Zhang04}), we have
    \begin{equation*}
    	\det G=(u'(s))^{n-1}(u'(s)+su''(s)).
    \end{equation*}
        The inverse matrix  $G^{-1}=(g^{j\overline{k}})$ can be obtained by computing the adjunct of $G$,
    \begin{equation*}
    	g^{j\overline{k}}=e^{-f}(u'(s))^{n-2}[(u'(s)+su''(s))\delta_{jk}-u''(s)\overline{z}_kz_j]
    \end{equation*}
    with $f=\log \det G$. Moreover we have
    \begin{equation*}
    	\partial \overline{\partial}f=\sum_{j,k=1}^{n}(\delta_{jk}f'(s)+f''(s)\overline{z}_jz_k)dz_j \wedge d\overline{z}_k
    \end{equation*}
    and the scalar curvature
    \begin{equation*}
    \begin{split}
    	R &=-\sum_{j,k=1}^{n} g^{j \overline{k}}\frac{\partial^2 f}{\partial z_j \partial \overline{z}_k} \\
    	  &=-e^{-f}(u'(s))^{n-2}(nu'(s)f'(s)+su'(s)f''(s)+(n-1)su''(s)f'(s)) \\
    	  &=-s^{1-n}e^{-f}[s^n(u'(s))^{n-1}f'(s)]'
    \end{split}
    \end{equation*}
    
    To summarize, $u(|z_1|^2+|z_2|^2+\cdot\cdot\cdot+|z_n|^2)$ is the potential of a K\"ahler metric on $\mathbb{C}^n\backslash\{0\}$ with scalar curvature $R$ if and only if the smooth function $u:(0,+\infty) \rightarrow \mathbb{R}$ satisfies the following system of ordinary differential equations:
    \begin{equation}\label{main equation}
    		\left\{
    	\begin{array}{ccc}
    	u'(s)>0,u'(s)+su''(s)>0 \\
    	e^{f(s)}=(u'(s))^{n-1}(u'(s)+su''(s))\\
    	s^{n-1}e^{f(s)}R=-[s^n(u'(s))^{n-1}f'(s)]'
    	\end{array}
    	\right.
    \end{equation}
   For later convenience we introduce an auxiliary function $g(s)=su'(s)$ and reformulate it as
   \begin{equation}\label{main equation 1}
   \left\{
   \begin{array}{ccc}
   g(s)>0,g'(s)>0 \\
   s^{n-1}e^{f(s)}=(g(s))^{n-1}g'(s)\\
   s^{n-1}e^{f(s)}R=-[s(g(s))^{n-1}f'(s)]'
   \end{array}
   \right.
   \end{equation}
   
  We want to emphasize that $g:(0,+\infty) \rightarrow \mathbb{R}$ satisfying the system of equation (\ref{main equation 1}) are positive and strictly increasing. It follows that $A=\lim\limits_{s \rightarrow 0+} g(s)$ and $B=\lim\limits_{s \rightarrow +\infty} g(s)$ always make sense. Moreover $A$ is finite and nonnegative, $B$ is positive finite or $+\infty$.
  
  If the scalar curvature R is constant, the system of ODE (\ref{main equation 1}) can be reduced to a single equation:
  \begin{lem}\label{potential equation}
  	If $u:(0,+\infty) \rightarrow \mathbb{R}$ is a smooth function such that $u(|z_1|^2+|z_2|^2+\cdots+|z_n|^2)$ is the potential of a K\"ahler metric with constant scalar curvature $R$ on $\mathbb{C}^n\backslash\{0\}$, then there exist constants $\lambda,\mu \in \mathbb{R}$ such that
  	\begin{equation}\label{pequation01}
  	\frac{g^{n-1}(s)g'(s)}{-\frac{R}{n(n+1)}g^{n+1}(s)+g^n(s)+\lambda g(s)+\mu}=\frac{1}{s}.
  	\end{equation}
  \end{lem}
  \pr Substitute the second equation in  (\ref{main equation 1}) into the third one we obtain:
  \begin{equation*}
  [sg^{n-1}(s)f'(s)]'=-\frac{R}{n}[g^n(s)]'
  \end{equation*}
  It follows that there exists a constant $\lambda$ such that
  \begin{equation}\label{be00}
  sg^{n-1}(s)f'(s)=-\frac{R}{n}g^n(s)+\lambda
  \end{equation}
  Differentiate the second equation in (\ref{main equation 1}) and then multiply both sides by $s$ we have:
  \begin{equation}\label{zero II}
  (n-1)s^{n-1}e^{f(s)}+s^{n}e^{f(s)}f'(s)=s[g^{n-1}(s)g'(s)]'
  \end{equation}
  Put in  (\ref{be00}) and the second equation of (\ref{main equation 1}), we obtain
  \begin{equation*}
  [g^n(s)]'-\frac{R}{n(n+1)}[g^{n+1}(s)]'+\lambda g'(s)=[sg^{n-1}(s)g'(s)]'
  \end{equation*}
  Hence there exists a constant $\mu$ such that
  \begin{equation*}
  \frac{g^{n-1}(s)g'(s)}{-\frac{R}{n(n+1)}g^{n+1}(s)+g^n(s)+\lambda g(s)+\mu}=\frac{1}{s}
  \end{equation*}         \qed
 
To end this section we oberve that smooth and singular K\"ahler metrics coming from the equation (\ref{main equation}) are distinguished by the following property:
\begin{lem}\label{smooth}
	If $u:[0,+\infty) \rightarrow \mathbb{R}$ is a smooth function such that $u(|z_1|^2+|z_2|^2+\cdot\cdot\cdot+|z_n|^2)$ is the potential of a K\"ahler metric on $\mathbb{C}^n$ and $g(s)=su'(s)$. Then
	\begin{itemize}
		\item[(1)] $\lim\limits_{s\rightarrow 0+} g(s)=0$;
		\item[(2)] $\lim\limits_{s \rightarrow 0+} sg'(s)=0$.
	\end{itemize} 
\end{lem}

\pr The potential $\varphi(z_1,z_2,...,z_n)=u(|z_1|^2+|z_2|^2+\cdot\cdot\cdot+|z_n|^2)$ is a smooth function on $\mathbb{C}^n$. It follows that $f(x)=\varphi(x,0,...,0)=u(x^2)$ is a smooth function on $\mathbb{R}$ and $u(s)=f(\sqrt{s})$. By a direct computation we obtain:
\[
g(s)=su'(s)=\frac{1}{2}\sqrt{s}f'(\sqrt{s}), \;sg'(s)=\frac{\sqrt{s}}{4}f'(\sqrt{s})+\frac{s}{4}f''(\sqrt{s}).
\] 
Then the lemma follows immediately.   \qed

\section{rotation invariant csck: smooth case}

In this section we solve the ODE (\ref{pequation01}) in the smooth case and give a proof of Theorem \ref{n smooth solution}.

\textit{Proof of Theorem \ref{n smooth solution}}: Suppose that $u:[0,+\infty) \rightarrow \mathbb{R}$ is a smooth function such that $u(|z|^2)$ is the potential of a K\"ahler metric with constant scalar curvature $R$ on $\mathbb{C}^n$. It follows from Lemma \ref{potential equation} that there exist constants $\lambda,\mu$ such that
\begin{equation*}
	-\frac{R}{n(n+1)}g^{n+1}(s)+g^n(s)+\lambda g(s)+\mu=sg^{n-1}(s)g'(s)
\end{equation*}
for $s \in (0,+\infty)$. Lemma \ref{smooth} implies that $\mu=0$ and $\lambda=0$. Then the equation (\ref{pequation01}) can be simplified as
\begin{equation}\label{smooth solution}
	\frac{g'(s)}{-\frac{R}{n(n+1)}g^2(s)+g(s)}=\frac{1}{s}
\end{equation}

\begin{itemize}
	\item[(1)] If $R=0$, it follows imediately that there exist constants $a,b$ with $a>0$ such that
	\begin{equation*}
		u(s)=as+b
	\end{equation*}
	and the K\"ahler form
	\begin{equation*}
		\omega=ai\sum_{j}^{n}dz_j \wedge d\overline{z}_j
	\end{equation*}
	
	\item[(2)] If $R=n(n+1)$, the equation (\ref{smooth solution}) can be rewritten as
	\[
	\frac{g'(s)}{-g(s)(g(s)-1)}=\frac{1}{s}
	\]
	or
	\[
		\frac{g'(s)}{g(s)}-\frac{g'(s)}{g(s)-1}=\frac{1}{s}
	\]
	It follows from Lemma \ref{no real roots} that $\lim\limits_{s \rightarrow 0+} g(s)=0$ and $\lim\limits_{s\rightarrow \infty} g(s)=1$ and there exists a constant $C_1$ such that
	\begin{equation*}
		\log g(s)-\log(1-g(s))=\log s-C_1
	\end{equation*}
	Then we have
	\begin{equation*}
		g(s)=\frac{s}{s+a}
	\end{equation*}
	with $a=e^{C_1}>0$. It follows that there exists a constant $c$ such that
	\begin{equation*}
		u(s)=\log(s+a)+c
	\end{equation*}
	and the K\"ahler form
	\begin{equation*}
		\omega=i\frac{\sum_{j=1}^{n}dz_j \wedge d \overline{z}_j}{\sum_{j=1}^{n}|z_j|^2+a}-i\frac{(\sum_{j=1}^{n}\overline{z}_jdz_j)\wedge \overline{(\sum_{j=1}^{n}\overline{z}_jdz_j)}}{(\sum_{j=1}^{n}|z_j|^2+a)^2}.
	\end{equation*}
	
	\item[(3)] We will show that the equation (\ref{smooth solution}) has no solution if $R=-n(n+1)$. 
	
	If $R=-n(n+1)$, the equation (\ref{smooth solution}) can be rewritten as
	\begin{equation*}
		\frac{g'(s)}{(g(s)+1)g(s)}=\frac{1}{s}
	\end{equation*}
It follows from Lemma \ref{no real roots} that it does not admit required solution.
\end{itemize}   \qed

\section{rotation invariant csck on $\mathbb{C}^2\backslash\{0\}$} 

    In this section we solve the ordinary differential equation (\ref{pequation01}) for dimension $n=2$ and give a proof of Theorem \ref{2 zero singular solution} and Theorem \ref{2 negative singular solution}. We also prove Theorem \ref{positive} partially.
    
\textit{Proof of Theorem \ref{2 zero singular solution}:}  
 It follows from Theorem \ref{potential equation} that there exist constants $\lambda,\mu$ such that
 \begin{equation}\label{2-zero}
 	\frac{g(s)g'(s)}{g^2(s)+\lambda g(s)+\mu}=\frac{1}{s}
 \end{equation}
 \begin{itemize}
 	\item[Case 1] $\mu=0, \lambda=0$.
 	
 	There exist constants $a,b$ with $a>0$ such that
 	\begin{equation*}
 		u(s)=as+b.
 	\end{equation*}
 	                   
 	\item[Case 2] $\mu = 0, \lambda \neq 0$.
 	
 	The equation (\ref{2-zero}) can be rewritten as 
 	               \begin{equation}
 	               	\frac{g'(s)}{g(s)+\lambda}=\frac{1}{s}
 	               \end{equation}
 	               
 	               	It follows from Lemma \ref{no real roots} that $A=-\lambda$. Then $\lambda<0$ and there exists a constant $C$ such that
 	               	\begin{equation*}
 	               	\log (g(s)+\lambda) =\log s +C
 	               	\end{equation*}
 	and
 	               	\begin{equation*}
 	               		g(s)=as+b.
 	               	\end{equation*}
	 	with $a=e^C>0, b=-\lambda>0$. It follows that there exists a constant $c$ such that
 	               	\begin{equation*}
 	               		u(s)=as+b\log s+c.
 	               	\end{equation*}

 	\item[Case 3] $\mu \neq 0$ and the polynomial $x^2+\lambda x+\mu$ has two distinct real roots $\alpha < \beta$.
 	
 	We have $\alpha\beta \neq 0$ and the equation (\ref{2-zero}) can be written as
 	\begin{equation*}
 			\frac{\beta}{\beta-\alpha}\frac{g'(s)}{g(s)-\beta} -\frac{\alpha}{\beta-\alpha}\frac{g'(s)}{g(s)-\alpha} =\frac{1}{s}
 	\end{equation*}
 	It follows from Lemma \ref{no real roots}  that $A=\beta$. Hence $\beta>0$ and there exists a constant $c$ such that
         	\begin{equation}\label{2-zero-(3)}
         	\frac{\beta}{\beta-\alpha}\log(g(s)-\beta)-\frac{\alpha}{\beta-\alpha}\log (g(s)-\alpha)=\log s+c
         	\end{equation}
	
         	On the other hand, Lemma \ref{existence} implies that there exists a unique smooth strictly increasing function $g(s)=su'(s)$ ranging from $\beta$ to $+\infty$ on $(0,+\infty)$ satisfying (\ref{2-zero-(3)}) if
	\[
	\alpha \neq 0, \beta>0, \alpha<\beta
	\]         
     
 	\item[Case 4] $\mu \neq 0$ and the polynomial $x^2+\lambda x+\mu$ has double real roots $\alpha$.
 	
 	We have $\alpha \neq 0$ and the equation (\ref{2-zero}) can be rewritten as
 	\begin{equation*}
 		\frac{g'(s)}{g(s)-\alpha}+\frac{\alpha g'(s)}{(g(s)-\alpha)^2}=\frac{1}{s}
 	\end{equation*}
         It follows from Lemma \ref{no real roots} that $A=\alpha$. Hence $\alpha>0$ and there exists a constant $c$ such that
 		\begin{equation}\label{2-zero-(4)}
 		\log(g(s)-\alpha)-\frac{\alpha}{g(s)-\alpha} =\log s+c
 		\end{equation}
 		
 		On the other hand, Lemma \ref{existence} implies that there exists a unique smooth strictly increasing function $g(s)=su'(s)$ ranging from $\alpha$ to $+\infty$ on $(0,+\infty)$ satisfying (\ref{2-zero-(4)}) if $\alpha>0$.

 	\item[Case 5] $\mu \neq 0$ and the polynomial $x^2+\lambda x+\mu$ has no real roots.
 	
 	      According to Lemma \ref{no real roots}, the equation (\ref{2-zero}) does not admit the required solution.
 \end{itemize}   \qed

   In general the potentials belong to the third and forth class above cannot be expressed by explicit formulas. But in some special cases, $\beta=-\alpha>0$ in the third class for example, one can solve that
   \begin{equation*}
   		u(s)=a(\sqrt{s^2+b^2}+b\log \frac{s}{\sqrt{s^2+b^2}+b})+c
   \end{equation*}
  where $a=e^C>0,b=\beta e^{_C}>0$ and $c$ is a constant.

\textit{Proof of Theorem \ref{2 negative singular solution}:}  It follows from Theorem \ref{potential equation} that there exist constnts $\lambda,\mu$ such that
\begin{equation}\label{2-negative}
	-\frac{g(s)g'(s)}{g^3(s)-g^2(s) -\lambda g(s)-\mu}=\frac{1}{s}
\end{equation}

It follows from Lemma \ref{no real roots} that $B=\lim\limits_{s \rightarrow +\infty} g(s)<+\infty$.
\begin{itemize}
	\item[Case 1] $\mu=0,\lambda=0$.
	
	The equation \ref{2-negative} can be rewritten as
	\begin{equation*}
		\frac{g'(s)}{g(s)}-\frac{g'(s)}{g(s)-1}=\frac{1}{s}
	\end{equation*}
By Lemma \ref{no real roots} we have $A=0$, then there exists a constant $C$ such that
		\begin{equation*}
			\log g(s)-\log (1-g(s))=\log s-C
		\end{equation*}
		and
		\begin{equation*}
			g(s)=\frac{s}{s+a}
		\end{equation*}
		with $a=e^C>0$. Hence there exists a constant $c$ such that
		\[
			u(s)=\log (s+a)+c.
		\]
		
	\item[Case 2] $\mu =0,\lambda \neq 0$ and the polynomial $x^2-x-\lambda$ has two distinct real roots $\alpha<\beta$.
	      
	    Then we have $\alpha \beta \neq 0,\alpha+\beta=1$ and the equation (\ref{2-negative}) can be rewritten as
	    \begin{equation*}
	    \frac{g'(s)}{g(s)-\alpha}-\frac{g'(s)}{g(s)-\beta}=\frac{\beta-\alpha}{s}
	    \end{equation*}
	
	    By Lemma \ref{no real roots} we have $A=\alpha, B=\beta$. Hence $\alpha>0$ and there exists a constant $\gamma$ such that
	    \begin{equation*}
	    \log(g(s)-\alpha)-\log(\beta-g(s))=(\beta-\alpha)\log s -\gamma
	    \end{equation*}
	    and
	    \begin{equation*}
	    g(s)=\frac{1}{2}(k+1)-\frac{k a}{s^{k}+a}
	    \end{equation*}
	    where $k=\beta-\alpha$ and $a=e^{\gamma}$ satisfying $a>0, 0<k<1$.

	\item[Case 3] $\mu =0,\lambda \neq 0$ and the polynomial $x^2-x-\lambda$ has double real roots or no real roots.
	
	  It follows from Lemma \ref{no real roots} (1) that the equation (\ref{2-negative}) does not admit the required solution.

\item[Case 4] $\mu \neq 0$ and the polynomial $x^3-x^2-\lambda x-\mu$ has three distinct real roots $\alpha<\beta<\gamma$.

Then $\alpha+\beta+\gamma=1,\alpha \beta \gamma \neq 0$ and the equation \ref{2-negative} can be rewritten as
\begin{equation*}
	\begin{split}
	-\frac{\alpha(\gamma-\beta)g'(s)}{g(s)-\alpha}+\frac{\beta(\gamma-\alpha)g'(s)}{g(s)-\beta}
	-\frac{\gamma(\beta-\alpha)g'(s)}{g(s)-\gamma}=\frac{(\beta-\alpha)(\gamma-\beta)(\gamma-\alpha)}{s}
	\end{split}
\end{equation*}
It follows from Lemma \ref{no real roots} that $A=\beta$ and $B=\gamma$. Hence $\beta >0$ and there exists a constant $c$ such that
	\begin{equation}\label{2-negative-(3)}
	\begin{split}
	&-\alpha(\gamma-\beta)\log(g(s)-\alpha)+\beta(\gamma-\alpha)\log(g(s)-\beta) \\
	&-\gamma(\beta-\alpha)\log(\gamma-g(s)) 
	=(\beta-\alpha)(\gamma-\beta)(\gamma-\alpha)\log s+c
	\end{split}	
	\end{equation}	
	On the other hand, according to Lemma \ref{existence} there exists a unique smooth strictly increasing function $g(s)=su'(s)$ ranging from $\beta$ to $\gamma $ on $(0,+\infty)$ satisfying (\ref{2-negative-(3)}) if 
	\[
	\alpha<\beta<\gamma, \alpha \neq 0, \beta >0, \alpha+\beta+\gamma=1.
	\]

\item[Case 5] $\mu \neq 0$ and $x^3-x^2-\lambda x-\mu$ has three real roots $\alpha,\beta,\beta$ with $\alpha \neq \beta$.

Then $\alpha \beta \neq 0,\alpha+2\beta=1$ and the equation (\ref{2-negative}) can be written as
\begin{equation*}
\frac{\alpha g'(s)}{g(s)-\beta}-\frac{\alpha g'(s)}{g(s)-\alpha}-\frac{\beta(\beta-\alpha)g'(s)}{(g(s)-\beta)^2}=\frac{(\beta-\alpha)^2}{s}
\end{equation*}
 It follows from Lemma \ref{no real roots} that $A=\beta, B=\alpha$. Hence $\alpha>\beta>0$ and  there exists a constant $c$ such that
	\begin{equation}\label{2-negative-(4)}
	\alpha\log(g(s)-\beta)-\alpha\log(\alpha-g(s))+\frac{\beta(\beta-\alpha)}{g(s)-\beta}=(\beta-\alpha)^2\log s+c.
	\end{equation}
By Lemma \ref{existence} there exists a unique smooth strictly increasing function $g(s)=su'(s)$ ranging from $\beta$ to $\alpha$ on $(0,+\infty)$ satisfying (\ref{2-negative-(4)}) if
\[
\alpha>\beta>0, \alpha+2\beta=1.
\]

\item[Case 6] $\mu \neq 0$ and $x^3-x^2-\lambda x-\mu$ has three real roots $\alpha,\alpha,\alpha$ or only one real root.

  It follows from lemma \ref{no real roots} that the equation (\ref{2-negative}) dose not admit the required solution.
	
\end{itemize}    \qed

\begin{thm}
	Let $u:(0,+\infty) \rightarrow \mathbb{R}$ be a smooth function such that $u(|z_1|^2+|z_2|^2)$ is the potential of a K\"ahler metric $\omega$ with constant scalar curvature $R$ on $\mathbb{C}^2-\{0\}$. Then $R$ cannot be $-6$.
\end{thm}
\pr If the K\"ahler metric $\omega$ has constant scalar curvature $-6$, it follows from Theorem \ref{potential equation} that there exist constants $\lambda,\mu$ such that

\begin{equation}\label{2-positive}
	\frac{g(s)g'(s)}{g^3(s)+g^2(s)+\lambda g(s)+\mu}=\frac{1}{s}
\end{equation}

It follows from Lemma \ref{no real roots} that $B=\lim\limits_{s \rightarrow +\infty} g(s)<+\infty$.
\begin{itemize}
	\item[Case 1] $\mu=0$.
	
	It follow from Lemma \ref{no real roots} that the polynomial $x^2+x+\lambda$ admits at least two nonnegative real roots. It is impossible.

	\item[Case 2] $\mu \neq 0$ and $x^3+x^2+\lambda x+\mu$ has three distinct real roots $\alpha<\beta<\gamma$.
	
	We have $\alpha+\beta+\gamma=-1,\alpha\beta\gamma \neq 0$. It follows from Lemma \ref{no real roots} that
	\[
	A=\alpha, B=\beta
	\]
	and
	\[
	0<\alpha < \beta < \gamma, \alpha+\beta+\gamma=-1
	\]
	It is impossible.

	\item[Case 3] $\mu \neq 0$ and $x^3+x^2+\lambda x+\mu$ has three real roots $\alpha,\beta,\beta$ with $\alpha \neq \beta$.
	
	We have $\alpha+2\beta=-1,\alpha \beta \neq 0$. It follows from Lemma \ref{no real roots} that
	\[
	A=\alpha, B=\beta
	\]
	and
	\[
	0<\alpha < \beta, \alpha+2\beta=-1
	\]	
        It is impossible.
        
	\item[Case 4] $\mu \neq 0$ and the polynomial $x^3+x^2+\lambda x+\mu$ has three real roots $\alpha,\alpha,\alpha$ or only one real root.
	
	 It follows from Lemma \ref{no real roots} that the equation (\ref{2-positive}) dose not admit the required solution.

\end{itemize}    
In conclusion, the constant scalar curvature of the K\"ahler metric $\omega=i\partial \overline{\partial}u$ can not be $-6$.
\qed

\begin{cor}\label{positive 2}
	There does not exist rotation invariant K\"ahler metric with  negative constant scalar curvature  on $\mathbb{C}^2-\{0\}$.
\end{cor}

\section{rotation invariant csck on $\mathbb{C}^3\backslash\{0\}$}

In this section we solve the equation (\ref{pequation01}) for $n=3$ and give proofs for Theorem \ref{3 zero singular solution} and Theorem \ref{3 negative singular solution}. We also complete the proof of Theorem \ref{positive}.

\textit{Proof of Theorem \ref{3 zero singular solution}:} It follows from theorem \ref{potential equation} that there exist constants $\lambda,\mu$ such that
\begin{equation}\label{3 zero scalar curvature}
	\frac{g^2(s)g'(s)}{g^3(s)+\lambda g(s)+\mu}=\frac{1}{s}
\end{equation}

\begin{itemize}
	\item[Case 1] $\mu=0,\lambda=0$.
	
	The equation (\ref{3 zero scalar curvature}) can be rewritten as
	\begin{equation*}
		\frac{g'(s)}{g(s)}=\frac{1}{s}
	\end{equation*}	
	and there exist constants $a,b$ with $a>0$ such that
	\begin{equation*}
		u(s)=as+b.
	\end{equation*}
	
	\item[Case 2] $\mu=0,\lambda \neq 0$.
	
	The equation (\ref{3 zero scalar curvature}) can be rewritten as
	\begin{equation*}
		\frac{g(s)g'(s)}{g^2(s)+\lambda}=\frac{1}{s}
	\end{equation*}
	It follows from Lemma \ref{no real roots} that $\lambda=-A^2<0$ with $A>0$ and $B=+\infty$. Then there exists a constant $C$ such that
		\begin{equation*}
			\frac{1}{2}\log (g^2(s)-A^2)=\log s+C
		\end{equation*}
	Hence
		\begin{equation*}
			g(s)=a\sqrt{s^2+b^2},
		\end{equation*}
	with $a=e^C>0, b=\frac{A}{a}>0$. It follows that there exists a constant $c$ such that
	\[
	u(s)=a(\sqrt{s^2+b^2}+b\log \frac{s}{\sqrt{s^2+b^2}+b})+c.
	\]	
	\item[Case 3] $\mu \neq 0$ and $x^3+\lambda x+\mu$ has three distinct real roots $\alpha < \beta < \gamma$.
	
	              We have 
	              \begin{equation*}
	              	\alpha+\beta+\gamma=0,\alpha\beta\gamma \neq 0,\alpha<0
	              \end{equation*}
	              It follows from Lemma \ref{no real roots} that $A=\gamma,B=+\infty$. The equation (\ref{3 zero scalar curvature}) can be rewritten as 
	              
	              \begin{equation*}
	              	\frac{\alpha^2(\gamma-\beta)g'(s)}{g(s)-\alpha}-\frac{\beta^2(\gamma-\alpha)g'(s)}{g(s)-\beta}+\frac{\gamma^2(\beta-\alpha)g'(s)}{g(s)-\gamma}=\frac{(\beta-\alpha)(\gamma-\beta)(\gamma-\alpha)}{s}
	              \end{equation*}
	       Then there exists a constant $c$ such that
	             \begin{equation}\label{3-(3)}
	             \begin{split}
	             	&\alpha^2(\gamma-\beta)\log(g(s)-\alpha)-\beta^2(\gamma-\alpha)\log(g(s)-\beta)  \\
		&+\gamma^2(\beta-\alpha)\log(g(s)-\gamma) 
		=(\beta-\alpha)(\gamma-\beta)(\gamma-\alpha)\log s+c
		\end{split}
	             \end{equation}
	             On the other hand, Lemma \ref{existence} implies that there exists a unique smooth strictly increasing function $g(s)$ ranging from $\gamma$ to $+\infty$ on $(0,+\infty)$  satisfying (\ref{3-(3)}) if
	             \[
	             \alpha<0,\beta \neq 0, \gamma>0,\alpha<\beta<\gamma, \alpha+\beta+\gamma =0.
	             \]
	             
	\item[Case 4] $\mu \neq 0$ and $x^3+\lambda x+\mu$ has three real roots $\alpha,\alpha,\beta$ with $\alpha \neq \beta$.

                 	Then $2\alpha+\beta=0,\alpha \beta \neq 0$ and the equation (\ref{3 zero scalar curvature}) can be rewritten as
                 \begin{equation*}
                 	\frac{5}{9}\frac{g'(s)}{g(s)-\alpha}+\frac{4}{9}\frac{g'(s)}{g(s)+2\alpha}+\frac{\alpha}{3}\frac{g'(s)}{(g(s)-\alpha)^2}=\frac{1}{s}
                 \end{equation*}
                 
                 Lemma \ref{no real roots} implies that $A=\alpha$ or $-2\alpha$. 
                 
                 If $A=-2\alpha$, then $\alpha<0,\beta=-2\alpha>0$. It follows from Lemma \ref{no real roots} that $B=+\infty$. Then there exists a constant $c$ such that
                 \begin{equation}\label{3-0-(4)}
                 	\frac{5}{9}\log (g(s)-\alpha)+\frac{4}{9}\log(g(s)+2\alpha)-\frac{\alpha}{3}\frac{1}{g(s)-\alpha}=\log s+c
                 \end{equation}
               On the other hand, Lemma \ref{existence} implies that there exists a unique smooth strictly increasing function $g(s)$ ranging from $-2\alpha$ to $+\infty$ on $(0,+\infty)$ satisfying (\ref{3-0-(4)}) if $\alpha<0$.

                 If $A=\alpha$, then $\alpha>0, \beta=-2\alpha<0$. It follows from Lemma \ref{no real roots} that $B=+\infty$ there exists a constant $c$ such that
                 \begin{equation}\label{3-(0)-5}
                 	\frac{5}{9}\log(g(s)-\alpha)+\frac{4}{9}\log (g(s)+2\alpha)-\frac{\alpha}{3}\frac{1}{g(s)-\alpha}=\log s+c
                 \end{equation}
         
              On the other hand, Lemma \ref{existence} implies that there exists a unique smooth strictly increasing function $g(s)$ ranging from $\alpha$ to $+\infty$ on $(0,+\infty)$ satisfying (\ref{3-(0)-5}) if $\alpha>0$.
                    
      \item[Case 5] $\mu \neq 0$ and the polynomial $x^3+\lambda x+\mu$ has only one real root.
      
                 We have $x^3+\lambda x+\mu=(x-\alpha)[(x-\beta)^2+\gamma^2]$ with 
                 \[\alpha \neq 0,\gamma>0,\alpha+2\beta=0.
                 \]
                 The equation (\ref{3 zero scalar curvature}) can be rewritten as
                 \begin{equation*}
                \frac{\alpha^2g'(s)}{g(s)-\alpha}+\frac{(\beta^2+\gamma^2-2\alpha\beta)(g(s)-\beta)g'(s)}{(g(s)-\beta)^2+\gamma^2}+\frac{(\beta^3+\beta\gamma^2-\alpha\beta^2+\alpha\gamma^2)g'(s)}{(g(s)-\beta)^2+\gamma^2}=\frac{(\alpha-\beta)^2+\gamma^2}{s}
                 \end{equation*}

                 It follows from Lemma \ref{no real roots} that $A = \alpha>0, B=+\infty$. Then there exists a constant $c$ such that
                 \begin{equation}\label{3-0-(6)}
                 \begin{split}
                 &\alpha^2\log(g(s)-\alpha)  
                 +\frac{(\beta^2+\gamma^2-2\alpha\beta)}{2}\log[(g(s)-\beta)^2+\gamma^2] \\ 
                 &+\frac{(\beta^3+\beta\gamma^2-\alpha\beta^2+\alpha\gamma^2)}{\gamma}\arctan \frac{g(s)-\beta}{\gamma}=[(\alpha-\beta)^2+\gamma^2]\log s+c
                 \end{split}	
                 \end{equation}
                On the other hand, Lemma \ref{existence} implies that  there exists a unique smooth strictly increasing function $g(s)$ ranging from $\alpha$ to $+\infty$ on $(0,+\infty)$ satisfying (\ref{3-0-(6)}) if
               \[
               \alpha>0,\gamma>0, \alpha+2\beta=0.
               \]
\end{itemize}    \qed

\textit{Proof of Theorem \ref{3 negative singular solution}:}  It follows from Theorem \ref{potential equation} that there exist constants $\lambda,\mu$ such that
\begin{equation}\label{3-negative}
	\frac{g^2(s)g'(s)}{-g^4(s)+g^3(s)+\lambda g(s)+\mu}=\frac{1}{s}
\end{equation}

It follows from Lemma \ref{no real roots} that $B<+\infty$.
\begin{itemize}
	\item[Case 1] $\mu=0,\lambda=0$.
	
	The equation \ref{3-negative} can be rewritten as
	\begin{equation*}
		\frac{g'(s)}{g(s)}-\frac{g'(s)}{g(s)-1}=\frac{1}{s}
	\end{equation*}
It follows from Lemma \ref{no real roots} that $A=0, B=1$. Then there exists a constant $C$ such that
	\begin{equation*}
		\log g(s)-\log (1-g(s))=\log s+C
	\end{equation*}
	and 
	\begin{equation*}
		g(s)=\frac{s}{s+a}
	\end{equation*}
	with $a=e^{-C}>0$. Hence there exists a constant $b$ such that
	\begin{equation*}
		u(s)=\log(s+a)+b.
	\end{equation*}
	
	\item[Case 2] $\mu=0,\lambda \neq 0$ and the polynomial $x^3-x^2-\lambda$ has three distinct real roots $\alpha<\beta<\gamma$.
	
	             We have 
	             \begin{equation*}
	             	\alpha+\beta+\gamma=1,\alpha\beta+\beta\gamma+\gamma\alpha=0,\alpha \beta \gamma \neq 0
	             \end{equation*} 
	             In particular $\alpha<0$ and the equation (\ref{3-negative}) can be rewritten as
	             \begin{equation*}
	             	-\frac{\alpha(\gamma-\beta)g'(s)}{g(s)-\alpha}+\frac{\beta(\gamma-\alpha)g'(s)}{g(s)-\beta}-\frac{\gamma(\beta-\alpha)g'(s)}{g(s)-\gamma}=\frac{(\beta-\alpha)(\gamma-\beta)(\gamma-\alpha)}{s}
	             \end{equation*}
	             
	             It follows from Lemma \ref{no real roots} that $A=\beta, B=\gamma$. Hence $\beta>0$ and there exists a constant $c$ such that
	       \begin{equation}\label{3-negative-(2)}
	       \begin{split}
	       &-\alpha(\gamma-\beta)\log(g(s)-\alpha)+\beta(\gamma-\alpha)\log(g(s)-\beta) \\
	      & -\gamma(\beta-\alpha)\log(\gamma-g(s))=(\beta-\alpha)(\gamma-\beta)(\gamma-\alpha)\log s+c
	       \end{split}
	       \end{equation}
	     On the other hand, Lemma \ref{existence} implies that there exists a unique smooth strictly increasing function $g(s)$ ranging from $\beta$ to $\gamma$ on $(0,+\infty)$ satisfying (\ref{3-negative-(2)}) if
	     \[
	     \alpha < 0<\beta < \gamma, \alpha+\beta+\gamma=1,\alpha\beta+\beta\gamma+\gamma\alpha=0.
	     \]

	\item[Case 3] $\mu=0,\lambda \neq 0$ and the polynomial $x^3-x^2-\lambda$ does not have  three distinct real roots.
	
	Then one of the following is true:
	\begin{enumerate}
	\item The polynomial  $x^3-x^2-\lambda$ has three real roots $\alpha, \alpha, \beta$ with $\alpha \neq \beta$. Then $\alpha=\frac{2}{3}, \beta =-\frac{1}{3}$;
	\item  The polynomial  $x^3-x^2-\lambda$ has only one  real root. 	
	\end{enumerate}
	It follows from Lemma \ref{no real roots} that the equation (\ref{3-negative}) does not admit the required solution.
	      	      
	      \item[Case 4] $\mu \neq 0$ and the polynomial $x^4-x^3-\lambda x-\mu$ has four distinct real roots $\alpha<\beta<\gamma<\delta$.
	      
	       We have 
	       \begin{equation*}
	       	\alpha\beta\gamma\delta \neq 0,\alpha+\beta+\gamma+\delta=1,\alpha\beta+\alpha\gamma+\alpha\delta+\beta\gamma+\beta\delta+\gamma\delta=0
	       \end{equation*}
	       The equation (\ref{3-negative}) can be rewritten as
	       \begin{equation*}
	       	\begin{split}	
		 &\frac{\alpha^2}{(\beta-\alpha)(\gamma-\alpha)(\delta-\alpha)}\frac{g'(s)}{g(s)-\alpha}-\frac{\beta^2}{(\beta-\alpha)(\gamma-\beta)(\delta-\beta)}\frac{g'(s)}{g(s)-\beta}   \\
	      &+\frac{\gamma^2}{(\gamma-\alpha)(\gamma-\beta)(\delta-\gamma)}\frac{g'(s)}{g(s)-\gamma}-\frac{\delta^2}{(\delta-\alpha)(\delta-\beta)(\delta-\gamma)}\frac{g'(s)}{g(s)-\delta} \\
	      &=\frac{1}{s}
	       \end{split}
	       \end{equation*}
	     Recall that we have $\alpha<0$. Then Lemma \ref{no real roots} implies that $A=\gamma>0, B=\delta$. Hence there exists a constant $c$ such that
	      \begin{equation}\label{3-negative-(3)}
	      \begin{split}
	      &\frac{\alpha^2}{(\beta-\alpha)(\gamma-\alpha)(\delta-\alpha)}\log(g(s)-\alpha)-\frac{\beta^2}{(\beta-\alpha)(\gamma-\beta)(\delta-\beta)}\log(g(s)-\beta)   \\
	      &+\frac{\gamma^2}{(\gamma-\alpha)(\gamma-\beta)(\delta-\gamma)}\log(g(s)-\gamma)-\frac{\delta^2}{(\delta-\alpha)(\delta-\beta)(\delta-\gamma)}\log(\delta-g(s)) \\
	      &=\log s+c
	      \end{split}
	      \end{equation}
	     On the other hand, Lemma \ref{existence} implies that there exists a unique smooth strictly increaing function $g(s)$ ranging from $\gamma$ to $\delta$ on $(0,+\infty)$ satisfying (\ref{3-negative-(3)}) if 
	     \[
	     \alpha< \beta < \gamma < \delta, \gamma >0, \alpha+\beta+\gamma+\delta=1, \alpha\beta+\alpha\gamma+\alpha\delta+\beta\gamma+\beta\delta+\gamma\delta=0, \alpha \beta\gamma\delta \neq 0.
	     \]

	      \item[Case 5] $\mu \neq 0$ and the polynomial $x^4-x^3-\lambda x-\mu$ has four real roots $\alpha,\alpha,\beta,\gamma$ with $\beta<\gamma$  and $\alpha,\beta,\gamma$ are distinct with each other.
	      
	      We have 
	      \begin{equation*}
	      	2\alpha+\beta +\gamma=1,\alpha\beta\gamma \neq 0,\alpha^2+2\alpha\beta+2\alpha\gamma+\beta\gamma=0
	      \end{equation*}
	      
	      and the equation (\ref{3-negative}) can be rewritten as
	      \begin{equation*}
	      	\begin{split}
	      	&-\frac{\gamma^2}{(\gamma-\beta)(\gamma-\alpha)^2}\frac{g'(s)}{g(s)-\gamma}+\frac{\beta^2}{(\gamma-\beta)(\beta-\alpha)^2}\frac{g'(s)}{g(s)-\beta} \\
	      	&-\frac{2\alpha\beta\gamma-\alpha^2(\beta+\gamma)}{(\gamma-\alpha)^2(\beta-\alpha)^2} \frac{g'(s)}{g(s)-\alpha}  
	      	-\frac{\alpha^2}{(\gamma-\alpha)(\beta-\alpha)}\frac{g'(s)}{(g(s)-\alpha)^2}=\frac{1}{s}
	      	\end{split}
	      \end{equation*}
	      
	     Recall that $\min \{\alpha,\beta,\gamma\}<0$, then Lemma \ref{no real roots} implies that the other two roots are positive and equal to $A, B$ respectively.      
	     
	     If $\alpha \notin \{A,B\}$, Lemma \ref{no real roots} implies that $A=\beta, B=\gamma$ and $\alpha<0<\beta<\gamma$. Then there exists a constant $c$ such that
	     
	     \begin{equation}\label{3-negative-(4)}
	     \begin{split}
	     &-\frac{\gamma^2}{(\gamma-\beta)(\gamma-\alpha)^2}\log(\gamma-g(s))+\frac{\beta^2}{(\gamma-\beta)(\beta-\alpha)^2}\log(g(s)-\beta) \\
	     &-\frac{2\alpha\beta\gamma-\alpha^2(\beta+\gamma)}{(\gamma-\alpha)^2(\beta-\alpha)^2} \log(g(s)-\alpha)  
	     +\frac{\alpha^2}{(\gamma-\alpha)(\beta-\alpha)}\frac{1}{g(s)-\alpha}=\log s+c
	     \end{split}
	     \end{equation}
	    On the other hand, Lemma \ref{existence} imply that there exists a unique smooth strictly increasing function $g(s)$ ranging from $\beta$ to $\gamma$ on $(0,+\infty)$ satisfying (\ref{3-negative-(4)}) If
	    \[
	    \alpha<0 < \beta<\gamma, 2\alpha+\beta+\gamma=1, \alpha^2+2\alpha\beta+2\alpha\gamma+\beta\gamma=0.
	    \]

	     If $\alpha \in \{A,B\}$, then $\beta<0$ and Lemma \ref{no real roots} imply that $A=\alpha, B=\gamma$ with $\beta<0<\alpha<\gamma$. Hence there exists a constant $c$ such that
	      \begin{equation}\label{3-negative-(5)}
	     \begin{split}
	     &-\frac{\gamma^2}{(\gamma-\beta)(\gamma-\alpha)^2}\log(\gamma-g(s))+\frac{\beta^2}{(\gamma-\beta)(\beta-\alpha)^2}\log(g(s)-\beta) \\
	     &-\frac{2\alpha\beta\gamma-\alpha^2(\beta+\gamma)}{(\gamma-\alpha)^2(\beta-\alpha)^2} \log(g(s)-\alpha)  
	     +\frac{\alpha^2}{(\gamma-\alpha)(\beta-\alpha)}\frac{1}{g(s)-\alpha}=\log s+c
	     \end{split}
	     \end{equation}

	     On the other hand, Lemma \ref{existence} implies that there exists a unique smooth strictly increasing function $g(s)$ ranging from $\alpha$ to $\gamma$ on $(0,+\infty)$ satisfying (\ref{3-negative-(5)}) if
	     \[
	     \beta<0<\alpha<\gamma, 2\alpha+\beta+\gamma=1, \alpha^2+2\alpha\beta+2\alpha\gamma+\beta\gamma=0.
	     \]

	      \item[Case 6] $\mu \neq 0$ and the polynomial $x^4-x^3-\lambda x-\mu$ has four real roots but at most two of them are distinct.
	         
	          Then one of the roots is negative and Lemma \ref{no real roots} imply that the equation (\ref{3-negative}) dose not admit the required solution.

	      \item[Case 7] $\mu \neq 0$ and the polynomial $x^4-x^3-\lambda x-\mu$ has two real roots $\alpha_1,\alpha_2$ with $\alpha_1 < \alpha_2$.
	      
	      Then $x^4-x^3-\lambda x-\mu=(x-\alpha_1)(x-\alpha_2)[(x-\beta)^2+\gamma^2]$ where
	      
	      \begin{equation*}
	      	\alpha_1+\alpha_2+2\beta=1,\beta^2+\gamma^2+2\alpha_1\beta+2\alpha_2\beta+\alpha_1\alpha_2=0,\gamma>0
	      \end{equation*} 
	      The equation (\ref{3-negative}) can be rewritten as
	      \begin{equation*}
	      	\begin{split}
	      	&\frac{1}{(\alpha_1-\beta)^2+\gamma^2}[\frac{\alpha^2_1 g'(s)}{g(s)-\alpha_1}+\frac{(\beta^2+\gamma^2-2\alpha_1\beta)(g(s)-\beta)g'(s)}{(g(s)-\beta)^2+\gamma^2}+\frac{[(\alpha_1+\beta)(\beta^2+\gamma^2)-2\alpha_1\beta^2]g'(s)}{(g(s)-\beta)^2+\gamma^2}]  \\
	      	&-\frac{1}{(\alpha_2-\beta)^2+\gamma^2}[\frac{\alpha^2_2 g'(s)}{g(s)-\alpha_2}+\frac{(\beta^2+\gamma^2-2\alpha_2\beta)(g(s)-\beta)g'(s)}{(g(s)-\beta)^2+\gamma^2}+\frac{[(\alpha_2+\beta)(\beta^2+\gamma^2)-2\alpha_2\beta^2]g'(s)}{(g(s)-\beta)^2+\gamma^2}]  \\
	      	&=\frac{\alpha_2-\alpha_1}{s}
	      	\end{split}
	      \end{equation*}
	      It follows from Lemma \ref{no real roots} that $A=\alpha_1, B=\alpha_2$ and $\alpha_2>\alpha_1>0$. Then there exists a constant $c$ such that
	      \begin{equation}\label{3-negative-(6)}
	      \begin{split}
	      &\frac{1}{(\alpha_1-\beta)^2+\gamma^2}[\alpha^2_1\log(g(s)-\alpha_1) +\frac{\beta^2+\gamma^2-2\alpha_1\beta}{2}\log((g(s)-\beta)^2+\gamma^2)  \\ &+\frac{(\alpha_1+\beta)(\beta^2+\gamma^2)-2\alpha_1\beta^2}{\gamma}\arctan\frac{g(s)-\beta}{\gamma}]  \\
	      &-\frac{1}{(\alpha_2-\beta)^2+\gamma^2}[\alpha^2_2\log(\alpha_2-g(s))+\frac{(\beta^2+\gamma^2-2\alpha_2\beta)}{2}\log((g(s)-\beta)^2+\gamma^2)   \\
	      &+\frac{(\alpha_2+\beta)(\beta^2+\gamma^2)-2\alpha_2\beta^2}{\gamma}\arctan\frac{g(s)-\beta}{\gamma}]  \\
	      &=(\alpha_2-\alpha_1)\log s+c
	      \end{split}
	      \end{equation}
	     On the other hand, Lemma \ref{existence} implies that there exists a unique smooth strictly increasing function $g(s)$ ranging from $\alpha_1$ to $\alpha_2$ on $(0,+\infty)$ satisfying (\ref{3-negative-(6)}) if
	     \[
	     0<\alpha_1<\alpha_2,\alpha_1+\alpha_2+2\beta=1,\beta^2+\gamma^2+2\alpha_1\beta+2\alpha_2\beta+\alpha_1\alpha_2=0,\gamma>0.
	     \]
	    
	      \item[Case 8] $\mu \neq 0$ and $x^4-x^3-\lambda x-\mu$ does not have two distinct real roots.
	      
	      It follows from Lemma \ref{no real roots} that the equation (\ref{3-negative}) dose not admit the required solution.
	      
\end{itemize}  \qed

It not hard to check that constants satisfying the conditions in $(2)-(5)$ in Theorem \ref{3 negative singular solution} exist. In fact we have the following data
\begin{itemize}
	\item[(2)] $\alpha=\frac{1-\sqrt{5}}{4},\beta=\frac{1}{2},\gamma=\frac{1+\sqrt{5}}{4}$;
	\item[(3)] $\alpha=\frac{1-\sqrt{21}}{8},\beta=\frac{1}{4},\gamma=\frac{1}{2},\delta=\frac{1+\sqrt{21}}{8}$;
	\item[(4)] $\alpha=-\frac{1}{6},\beta=\frac{1}{2},\gamma=\frac{5}{6}$;
	\item[(5)] $\alpha=\frac{1}{4},\beta=\frac{1-\sqrt{6}}{4},\gamma=\frac{1+\sqrt{6}}{4}$;
	\item[(6)] $\alpha_1=\frac{1}{2},\alpha_2=\frac{3}{2},\beta=-1,\gamma=\frac{3}{2}$
\end{itemize} 
satisfying the correspongding conditions.

\begin{thm}
	If $u(s):(0,+\infty) \rightarrow \mathbb{R}$ is a smooth function such that $u(|z_1|^2+|z_2|^2+|z_3|^2)$ is the potential of a K\"ahler metric $\omega$ with constant scalar curvature $R$ on $\mathbb{C}^3-\{0\}$. Then $R$ cannot be $-12$.
	
\end{thm}
\pr If $R=-12$, it follows from Theorem \ref{potential equation} that there exist constants $\lambda,\mu$ such that
     \begin{equation}\label{3-positive}
     	\frac{g^2(s)g'(s)}{g^4(s)+g^3(s)+\lambda g(s)+\mu}=\frac{1}{s}
     \end{equation}
  It follows from Lemma \ref{no real roots} that $B<+\infty$.
     \begin{itemize}
     	\item[Case 1] $\mu=0,\lambda=0$.
	
	It follows from Lemma \ref{no real roots} that the equation (\ref{3-positive}) does not admit the required solution.
     	
     	\item[Case 2] $\mu=0,\lambda \neq 0$.
	
	Then one of the following will be true:
	\begin{enumerate}
	\item The polynomial $x^3+x^2+\lambda$ has three distinct real roots $\alpha<\beta<\gamma$. In this case one has $\alpha<0$.
	\item The polynomial $x^3+x^2+\lambda$ has three real roots but at most two distinct ones. In this case the polynomial $x^3+x^2+\lambda$ does not admit two distinct positive real roots. 
	\item The polynomial $x^3+x^2+\lambda$ has only one real root.
		\end{enumerate}
	It follows from Lemma \ref{no real roots} that the equation (\ref{3-positive}) does not admit the required solution.
	         
        \item[Case 3] $\mu \neq 0$ and $x^4+x^3+\lambda x+\mu$ has four distinct real roots $\alpha<\beta <\gamma <\delta$.
        
         We have 
         \begin{equation*}
         	\alpha+\beta+\gamma+\delta=-1,\alpha\beta+\alpha\gamma+\alpha\delta+\beta\gamma+\beta\delta+\gamma\delta=0,\alpha\beta\gamma\delta \neq 0
         \end{equation*}
         Then we have $\alpha<0$. It follows from Lemma \ref{no real roots} that $A=\beta, B=\gamma$ and $\beta>0$. But Lemma \ref{3-positive-(1)} implies that there does not exist constants $\alpha,\beta,\gamma,\delta$ such that
         \[
         \alpha<0<\beta<\gamma<\delta,\alpha+\beta+\gamma+\delta=-1,\alpha\beta+\alpha\gamma+\alpha\delta+\beta\gamma+\beta\delta+\gamma\delta=0.
          \]
         
        \item[Case 4] $\mu \neq 0$ and $x^4+x^3+\lambda x+\mu$ has four  real roots $\alpha,\alpha,\beta ,\gamma $ with $\beta<\gamma$ and $\alpha,\beta,\gamma$ are distinct with each other.
        
        We have 
        \begin{equation*}
        	2\alpha+\beta+\gamma=-1,\alpha^2+2\alpha\beta+2\alpha\gamma+\beta\gamma=0,\alpha\beta\gamma \neq 0
        \end{equation*}
       Hence $\min \{\alpha,\beta,\gamma \}<0$ and Lemma \ref{no real roots} implies that the left two roots are positive and equal to $A, B$ respectively. 
       
       If $\alpha \notin \{A,B\}$, then  Lemma \ref{no real roots} implies that $A=\beta, B=\gamma$. It follows that $\alpha<0 <\beta<\gamma$ and the equation (\ref{3-positive}) dose not admit the required solution by Lemma \ref{no real roots} (1).
          
       If $\alpha \in \{A,B\}$, then Lemma \ref{no real roots} implies that $\beta<0 <\gamma<\alpha$ and $A=\gamma, B=\alpha$. But  Lemma \ref{3-positive-(2)} implies that constants $\alpha,\beta,\gamma$ satisfying all theses constraints do not exist.

        \item[Case 5] $\mu \neq 0$ and $x^4+x^3+\lambda x+\mu$ has four real roots but at most two of them are distinct.        
            
          Then at least one of these roots is negative and Lemma \ref{no real roots} implies that the equation (\ref{3-positive}) dose not admit the required solution.

        \item[Case 6] $\mu \neq 0$ and the polynomial $x^4+x^3+\lambda x+\mu$ has at most two real roots.
        
        Then one of the following is true:
        \begin{enumerate}
        \item The polynomial $x^4+x^3+\lambda x+\mu$ has two distinct real roots; In this case the polynomial is negative between theses two real roots.
        \item  The polynomial $x^4+x^3+\lambda x+\mu$ has two real roots $\alpha,\alpha$;
        \item  The polynomial $x^4+x^3+\lambda x+\mu$ dose not have real roots.
        \end{enumerate}
        It follows from Lemma \ref{no real roots} that the equation (\ref{3-positive}) does not admit the required solution.

  \end{itemize}\qed

\begin{cor}\label{positive 3}
	There dose not exist rotation invariant K\"ahler metric with negative constant scalar curvature on $\mathbb{C}^3\backslash\{0\}$.
\end{cor}
 
Theorem \ref{positive} follows from Corolary \ref{positive 2} and Corolary \ref{positive 3}.

\section{Technical lemmas}

This section is devoted to some technical Lemmas. Lemma \ref{existence} ensure the existence of those solutions obtained by implicit functions. A prior observation in Lemma \ref{no real roots} simplifies arguments to solve the ordinary differential equation (\ref{pequation01}) a lot. Lemma \ref{3-positive-(1)} and Lemma \ref{3-positive-(2)} are used to get rid of two difficult possibilities occurred in the solution for positive scalar curvature case in dimension three.

\begin{lem}\label{existence}
Let $h:(\alpha,\beta) \rightarrow \mathbb{R}$ be a smooth strictly increasing function with 
\[
\lim_{t \rightarrow \alpha} h(t)=-\infty, \lim_{t \rightarrow \beta} h(t)=+\infty
\]
Then for any constants $A>0$ and $C$, there exists a unique smooth strictly increasing function $g:(0,+\infty) \rightarrow \mathbb{R}$ such that
\[
h(g(s))=A\log s+C
\]
and $\lim\limits_{s \rightarrow 0+} g(s)=\alpha,\lim\limits_{s \rightarrow +\infty} g(s)=\beta$.
\end{lem}
\pr Obvious. \qed

\begin{lem}\label{no real roots}
If $H(x)$ is a polynomial of degree $m$ and the ordinary differential equation
\begin{equation}\label{no solution}
\frac{g^k(s)g'(s)}{H(g(s))} =\frac{1}{s}
\end{equation}
admits  a smooth solution $g(s)$ on $(0,+\infty)$ with $g(s)>0, g'(s)>0$. Denote by $A=\lim\limits_{s \rightarrow 0+} g(s), B=\lim\limits_{s \rightarrow +\infty} g(s)$. Then
\begin{enumerate}
\item $H(x)>0$ for $x \in (A,B)$ and $H(A)=0$.
\item If $B=+\infty$, then $\deg H \leq k+1$. Moreover $A$ is the largest and nonnegative real root of $H(x)$ and $H(x)$ is positive on $(A,+\infty)$.
\item If $B<+\infty$, then $H(B)=0$. Moreover $A,B$ are two succesive nonnegative real roots of the plolymomial $H(x)$ and $H(x)$ is positive on the interval $(A,B)$.
\end{enumerate}
\end{lem}
\pr Integrate the equation (\ref{no solution}) on the interval $s \in [s_1,s_2]$ we have:
\begin{equation}\label{appendix1}
\int_{g(s_1)}^{g(s_2)}\frac{x^k}{H(x)}dx=\log s_2-\log s_1
\end{equation}
for any $0<s_1< s_2$. 
\begin{enumerate}
\item  For any $x \in (A,B)$ we can choose $s_0 \in (0,+\infty)$ such that $x=g(s_0)$. Then $H(x)=H(g(s_0)))=s_0g^k(s_0)g'(s_0)>0$.

Suppose that $H(A) \neq 0$, then there exists a constant $C_1>0$ such that
\[
H(x) \geq \frac{1}{C_1}
\]
for all $x \in [A,g(1)]$. Hence the integral
\[
\int_{g(s_1)}^{g(1)} \frac{x^k}{H(x)}dx \leq \frac{C_1}{k+1}g^{k+1}(1)
\]
is uniformly bounded for $0<s_1<1$. It contradicts with the identity (\ref{appendix1}). It follows that $A$ is a nonnegative real root of $H(x)$.

\item  It follows from (1) and the fact $B=+\infty$ that $A$ is the largest nonnegative real root of $H(x)$. Then there exists a constant $C_2>0$ such that
\[
\frac{x^m}{H(x)} \leq C_2
\]
for $x \in [g(1),+\infty)$. 

If $\deg H > k+1$, the integral
\[
\int_{g(1)}^{g(s_2)} \frac{x^k}{H(x)}dx \leq C_2 \int_{g(1)}^{g(s_2)} \frac{1}{x^{m-k}} dx \leq \frac{C_2}{m-k-1}g^{k+1-m}(1)
\]
is uniformly bounded for $s_2>1$. It contradicts with the identity (\ref{appendix1}).

\item Suppose that $H(B) \neq 0$. Then there exists a constant $C_3>0$ such that
\[
H(x) \geq \frac{1}{C_3}
\]
for $x \in [g(1),B]$. Then the integral
\[
\int_{g(1)}^{g(s_2)}\frac{x^k}{H(x)}dx \leq C_3 \int_{g(1)}^{g(s_2)} x^k dx \leq \frac{C_3}{k+1}[B^{k+1}-g^{k+1}(1)]
\]
is uniformly bounded for $s_2>1$. It contradicts with the identity (\ref{appendix1}).
\end{enumerate}
\qed

\begin{lem}\label{3-positive-(1)}
	There does not exist real numbers $\alpha,\beta,\gamma,\delta$ such that
	\begin{itemize}
		\item[(1)] $\alpha+\beta+\gamma+\delta=-1$;
		\item[(2)] $\alpha\beta+\alpha\gamma+\alpha\delta+\beta\gamma+\beta\delta+\gamma\delta=0$;
		\item[(3)] $\alpha<0<\beta<\gamma<\delta$.
	\end{itemize}
\end{lem}
\pr We will prove that for real numbers $\alpha,\beta,\gamma,\delta$ satisfying (1),(3), we have
\begin{equation*}
	J(\alpha,\beta,\gamma,\delta)=\alpha\beta+\alpha\gamma+\alpha\delta+\beta\gamma+\beta\delta+\gamma\delta<0
\end{equation*}
In fact we have
\begin{equation*}
	\begin{split}
		J(\alpha,\beta,\gamma,\delta)&=\alpha(\beta+\gamma+\delta) +(\beta+\gamma)\delta+\beta\gamma  \\
		                             &=\alpha(-1-\alpha)+(\beta+\gamma)(-1-\alpha-\beta-\gamma)+\beta\gamma \\
		                             &=-\alpha^2-(1+\beta+\gamma)\alpha-(1+\beta+\gamma)(\beta+\gamma)+\beta\gamma
	\end{split}
\end{equation*}
It follows from $\gamma<\delta$ and (1) that $\alpha<-1-\beta-2\gamma <-\frac{1+\beta+\gamma}{2}$. Hence
\begin{equation*}
\begin{split}
J(\alpha,\beta,\gamma,\delta)&<-(-1-\beta-2\gamma)^2-(1+\beta+\gamma)(-1-\beta-2\gamma)-(1+\beta+\gamma)(\beta+\gamma)+\beta\gamma  \\
                             &=-\gamma(1+\beta+2\gamma)-(1+\beta+\gamma)(\beta+\gamma)+\beta\gamma  \\
                             &<0
\end{split}	
\end{equation*}   \qed

\begin{lem}\label{3-positive-(2)}
	There does not exist real numbers $\alpha,\beta,\gamma$ such that
	\begin{itemize}
		\item[(1)] $2\alpha+\beta+\gamma=-1$;
		\item[(2)] $\alpha^2+2\alpha\beta+2\alpha\gamma+\beta\gamma=0$;
		\item[(3)] $\beta<0<\gamma<\alpha$.
	\end{itemize}
\end{lem}
\pr We will prove that for real numbers $\alpha,\beta,\gamma$ satisfying $(1)(3)$, we have
\begin{equation*}
	I(\alpha,\beta,\gamma)=\alpha^2+2\alpha\beta+2\alpha\gamma+\beta\gamma<0
\end{equation*}

In fact we have
\begin{equation*}
	\begin{split}
	I(\alpha,\beta,\gamma)&=\beta(2\alpha+\gamma)+\alpha^2+2\alpha\gamma  \\
	                      &=\beta(-1-\beta)+ \frac{(-1-\beta-\gamma)^2}{4}+(-1-\beta-\gamma)\gamma\\
	                      &=-\frac{3}{4}\beta^2-\frac{1}{2}(\gamma+1)\beta+\frac{1}{4}(\gamma+1)^2-\gamma(\gamma+1)
	\end{split}
\end{equation*}
It follows from (1) and $\gamma<\alpha$ that $\beta<-1-3\gamma<-\frac{1}{3}(\gamma+1)$. Hence
\begin{equation*}
\begin{split}
I(\alpha,\beta,\gamma)&<-\frac{3}{4}(-1-3\gamma)^2-\frac{1}{2}(\gamma+1)(-1-3\gamma)+\frac{1}{4}(\gamma+1)^2-\gamma(\gamma+1) \\
                      &=-6\gamma^2-3\gamma \\
                      &<0
\end{split}
\end{equation*}  \qed

\end{document}